\documentclass{amsart}
\usepackage{times, pdfsync}

\usepackage{xypic}
\input xy
\xyoption{all}
\usepackage{epsfig}
\usepackage{amsthm}
\usepackage{amssymb}
\usepackage{amsmath}
\usepackage{amscd}
\usepackage{color}

\newcommand{\bg}{\begin{equation}}
\newcommand{\ed}{\end{equation}}
\newcommand{\bga}{\begin{eqnarray}}
\newcommand{\eda}{\end{eqnarray}}
\newcommand{\pf}{\textbf{Proof:\ }}

\def\cbdu{\par{\raggedleft$\Box$\par}}

\newtheorem {Theorem}  {Theorem}

\numberwithin{Theorem}{section}

\newtheorem {Lemma}[Theorem]  {Lemma}

\theoremstyle{definition}

\theoremstyle{remark}

\newtheorem {Corollary}[Theorem]{\bf Corollary}

\expandafter\chardef\csname pre amssym.def
at\endcsname=\the\catcode`\@ \catcode`\@=11
\def\undefine#1{\let#1\undefined}
\def\newsymbol#1#2#3#4#5{\let\next@\relax
 \ifnum#2=\@ne\let\next@\msafam@\else
 \ifnum#2=\tw@\let\next@\msbfam@\fi\fi
 \mathchardef#1="#3\next@#4#5}
\def\mathhexbox@#1#2#3{\relax
 \ifmmode\mathpalette{}{\m@th\mathchar"#1#2#3}%
 \else\leavevmode\hbox{$\m@th\mathchar"#1#2#3$}\fi}
\def\hexnumber@#1{\ifcase#1 0\or 1\or 2\or 3\or 4\or 5\or 6\or 7\or 8\or
 9\or A\or B\or C\or D\or E\or F\fi}

\font\teneufm=eufm10 \font\seveneufm=eufm7 \font\fiveeufm=eufm5
\newfam\eufmfam
\textfont\eufmfam=\teneufm \scriptfont\eufmfam=\seveneufm
\scriptscriptfont\eufmfam=\fiveeufm

\catcode`\@=\csname pre amssym.def at\endcsname

\newcounter{remark}
\renewcommand{\div}{\mbox{div}}

\def  \12  {{\frac{1}{2}}}

\def\build#1_#2^#3{\mathrel{\mathop{\kern 0pt#1}\limits_{#2}^{#3}}}

\begin{document}

\title[Norm inflation in Besov Spaces ]{Norm Inflation for Generalized Navier-Stokes Equations}


\author [Alexey Cheskidov]{Alexey Cheskidov}
\address{Department of Mathematics, Stat. and Comp.Sci.,  University of Illinois Chicago, Chicago, IL 60607,USA}
\email{acheskid@math.uic.edu} 
\author [Mimi Dai]{Mimi Dai}
\address{Department of Applied Mathematics, University of Colorado Boulder, Boulder, CO 80303,USA}
\email{mimi.dai@colorado.edu} 

\thanks{The work of Alexey Cheskidov was partially supported by NSF Grant
DMS--1108864.}





\begin{abstract}
We consider the incompressible Navier-Stokes equation with a fractional power $\alpha\in[1,\infty)$ of the Laplacian in the three dimensional case. We prove
the existence of a smooth solution with arbitrarily small in $\dot{B}_{\infty,p}^{-\alpha}$ ($2<p \leq \infty$) initial  data that becomes arbitrarily large in
$\dot{B}_{\infty,\infty}^{-s}$ for all $s> 0$ in arbitrarily small time. 
This extends the result of Bourgain and Pavlovi\'{c} \cite{BP}
for the classical Navier-Stokes equation which utilizes the fact that the energy transfer to low modes increases norms with negative smoothness indexes.
It is remarkable that the space $\dot{B}_{\infty,\infty}^{-\alpha}$ is supercritical
for $\alpha >1$. Moreover, the norm inflation occurs even in the case
$\alpha \geq 5/4$ where the global regularity is known.

\bigskip


\hspace{0.02cm}CLASSIFICATION CODE: 76D03, 35Q35.
\end{abstract}

\maketitle

\section{Introduction}

In this paper we study the three dimensional incompressible Navier-Stokes equations
with a fractional power of the Laplacian:
\begin{equation}\label{FNSE}
\begin{split}
u_t +(u\cdot\nabla) u+\nabla p=-\nu(-\Delta)^\alpha u,\\
\nabla \cdot u=0,\\
u(x, 0)=u_0, 
\end{split}
\end{equation}
where $x\in\mathbb{R}^3$, $t\geq 0$, $u$ is the fluid velocity, $p$
is the pressure of the fluid and $\nu>0$ is the kinematic viscosity coefficient. The initial data $u_0$ is divergence free.  The power $\alpha=1$ corresponds to the classical Navier-Stokes equations.
A vast amount of literature has been devoted to these equations, for background we refer the readers to \cite{CF} and \cite{Te}.

Solutions to the fractional Navier-Stokes equation (\ref{FNSE}) have the following scaling property.  If $(u(x, t), p(x,t))$ solves system (\ref{FNSE})
with the initial data $u_0(x)$, then
\begin{equation}
u_\lambda(x,t) =\lambda^{2\alpha-1} u(\lambda x, \lambda^{2\alpha}t), \ p_\lambda(x,t) =\lambda^{2(2\alpha-1)}
p(\lambda x, \lambda^{2\alpha}t) \notag
\end{equation}
solves the system (\ref{FNSE}) with the initial data
\begin{equation}
u_{0\lambda} =\lambda^{2\alpha-1} u_0(\lambda x). \notag
\end{equation}
A space that is
invariant under the above scaling is called a critical space.
The largest critical space in three dimension for the fractional NSE (\ref{FNSE}) is the Besov space $\dot B_{\infty,\infty}^{1-2\alpha}$ (see \cite{Ca}).

The study of the Navier-Stokes equations in
critical spaces has been a focus of the research activity
since the initial work of Kato \cite{Kato}.  In 2001, Koch and Tataru \cite{KT}
established the global well-posedness of the classical Navier-Stokes equations with
small initial data in the space $BMO^{-1}$. Then the question whether this result can be
extended to the largest critical space $\dot{B}^{-1}_{\infty,\infty}$ had become of great interest among
researchers, but it still remains open.

The first indication that such an extension might not be possible came in the work by
Bourgain and Pavlovi\'{c} \cite{BP} who showed the
norm inflation for the classical Navier-Stokes equations in $\dot{B}^{-1}_{\infty,\infty}$.
More precisely, they constructed arbitrarily small initial
data in $\dot{B}^{-1}_{\infty,\infty}$, such that mild solutions with this data become arbitrarily large in $\dot{B}^{-1}_{\infty,\infty}$ after an arbitrarily short time. This result was later extended to generalized Besov spaces smaller than $B^{-1}_{\infty,p}$, $p>2$ by Yoneda \cite{Yo}.
Moreover, in \cite{CS} Cheskidov and Shvydkoy proved the
existence of discontinuous Leray-Hopf solutions of the Navier-Stokes equations in  $\dot{B}^{-1}_{\infty,\infty}$ with arbitrarily small initial data.

Recently, Yu and Zhai \cite{YZ} considered the fractional Navier-Stokes equations \eqref{FNSE}
with $\alpha \in (1/2,1)$ and
showed global well-posedness for small initial data in the largest critical space
$\dot B_{\infty,\infty}^{1-2\alpha}$, conjecturing that the above mentioned ill-posedness
results could not be extended to the hypodissipative case $\alpha <1$.

Indeed, in the recent work \cite{CSnew} Cheskidov and Shvydkoy were able to prove
the existence of discontinuous Leray-Hopf solutions in the largest critical space 
with arbitrarily small initial data for  $\alpha\in[1,5/4)$. However, the construction broke down
for $\alpha <1$.

In this paper we consider the case $\alpha \in [1,\infty)$ and demonstrate
that the natural space for norm inflation is not critical, but
$\dot{B}_{\infty,p}^{-\alpha}$. Note that it is critical only in the classical case
$\alpha =1$, and it is not scaling invariant otherwise. More precisely, we prove
the existence of a smooth space-periodic solution with arbitrarily small in $\dot{B}_{\infty,p}^{-\alpha}$ ($2<p\leq \infty$) initial  data that becomes
arbitrarily large in $\dot{B}_{\infty,\infty}^{-s}$ for all $s> 0$ in arbitrarily small time. 
This recovers
Bourgain and Pavlovi\'{c}'s ill-posedness result in the case $\alpha =1$, and shows
the norm inflation in spaces $\dot{B}_{\infty,p}^{-s}$ for all $s\geq \alpha$, $p\in(2,\infty]$ in the case $\alpha\geq1$.
The case $\alpha >1$ is particularly interesting since the norm inflation occurs not only in critical spaces, but
even in supercritical spaces, which suggests that a small initial data result might be out of reach there.
It is remarkable that the norm inflation holds
even in the case $\alpha\geq 5/4$, where the global regularity is known.
In that case the smooth solution that exhibits the norm inflation can be extended globally in time.

Our construction is similar to the one of Bourgain and Pavlovi\'{c}, but we have
to deal with the lack of continuity of the bilinear operator corresponding to the fractional heat kernel on a modified Koch-Tataru adapted space.
This result is also based on the fact that a backwards energy cascade, harmless as far as the regularity of a solution is concerned, results in the
growth of Besov norms with negative smoothness indexes. In this construction we also make sure
that the initial data is space-periodic and has a finite energy if viewed on a torus.
Namely, we show that
\begin{Theorem}\label{Mthm}
Let $\alpha \geq 1$. For any $\delta >0$ and $2<p \leq \infty$  there exists
a smooth space-periodic solution $u(t)$ of (\ref{FNSE}) with period $2\pi$ and
the initial data
\begin{equation}\notag
\|u(0)\|_{\dot{B}^{-\alpha}_{\infty, p}} \lesssim \delta 
\end{equation}
that satisfies, for some $0<T<\delta$ and all $s> 0$,
\begin{equation}\notag
\|u(T)\|_{\dot{B}^{-s}_{\infty, \infty}}\gtrsim \frac{1}{\delta}.
\end{equation}
\end{Theorem}
We refer the reader to the  beginning of section of Preliminaries for the definition of the symbol $\lesssim$.

Note that the homogeneous and non-homogeneous Besov norms are equivalent
for periodic functions.  Therefore, for the space-periodic solution in Theorem~\ref{Mthm} we have 
\[
\|u(0)\|_{\dot{B}^{-s}_{\infty, p}} \lesssim \|u(0)\|_{B^{-s}_{\infty, p}}  \lesssim
\|u(0)\|_{\dot{B}^{-\alpha}_{\infty, p}} \qquad \text{for all} \qquad s \geq \alpha, 
\]  
Also, since
\[
\|u(T)\|_{\dot{B}^{-s}_{\infty, p}} \gtrsim \|u(T)\|_{\dot{B}^{-s}_{\infty, \infty}},
\]
the norm inflation occurs in all the spaces $\dot{B}^{-s}_{\infty, p}$,
$s \in [\alpha, \infty)$, $p\in (2, \infty]$. More precisely, we have the following.
\begin{Corollary}
Let $\alpha \geq 1$. For any  $s \geq \alpha$, $p\in(2, \infty]$, and $\delta >0$ there exists
a smooth space-periodic solution $u(t)$ of (\ref{FNSE}) with
the initial data
\begin{equation}\notag
\|u(0)\|_{\dot{B}^{-s}_{\infty, p}} \lesssim \delta, 
\end{equation}
that satisfies, for some $0<T<\delta$,
\begin{equation}\notag
\|u(T)\|_{\dot{B}^{-s}_{\infty, p}}\gtrsim \frac{1}{\delta}.
\end{equation}
\end{Corollary}
Moreover, due to the embedding of the Triebel-Lizorkin space $\dot{F}^{-s}_{\infty, p}$
\[
\dot{B}^{-s}_{\infty,p} \subset \dot{F}^{-s}_{\infty, p} \subset \dot{F}^{-s}_{\infty, \infty} =\dot{B}^{-s}_{\infty, \infty}, 
\]
Theorem~\ref{Mthm} also gives norm inflation in Triebel-Lizorkin spaces $\dot{F}^{-s}_{\infty, p}$, $s\geq \alpha$, $p \in (2,\infty]$:
\begin{Corollary}
Let $\alpha \geq 1$. For any  $s \geq \alpha$, $p\in(2, \infty]$, and $\delta >0$ there exists
a smooth space-periodic solution $u(t)$ of (\ref{FNSE}) with
the initial data
\begin{equation}\notag
\|u(0)\|_{\dot{F}^{-s}_{\infty, p}} \lesssim \delta, 
\end{equation}
that satisfies, for some $0<T<\delta$,
\begin{equation}\notag
\|u(T)\|_{\dot{F}^{-s}_{\infty, p}}\gtrsim \frac{1}{\delta}.
\end{equation}
\end{Corollary}

We now recall some  auxiliary concepts  related to the plane waves, which are necessary in the sequel:
\begin{itemize}
 \item  The ``diffusion" of a plane wave $v\cos (k\cdot x)$ in $R^3$ under the fractional Laplacian $-(-\Delta)^\alpha$ is given by
\[e^{-t(-\Delta)^\alpha }v\cos (k\cdot x)=e^{-|k|^{2\alpha} t}v\cos (k\cdot x)\] Thus the magnitude of the diffusion of a
plane wave dies down in time in the scale that is measured by $|k|^{2\alpha}$.
\item It is easy to  see
that $u = e^{-|k|^{2\alpha} t}v\cos (k\cdot x)$ solves the system (\ref{FNSE}) when the
wave vector $k$ is orthogonal to the amplitude vector $v$.
\item The
nonlinear interaction of two such diffusions in the system (\ref{FNSE}) can be
captured, and it produces only a slower diffusion if the two wave
vectors are close.
\end{itemize}
We note that these observations are the basis of the original argument  of Bourgain and
Pavlovi\'{c}  in \cite{BP}. We will use them to  construct a
combination of such ``diffusions"  with least nonlinear interactions yet
producing enough slower  ``diffusions"  to cause the norm inflation in
short time.

The rest of the paper is organized as: in Section \ref{sec:pre} we introduce some notations that shall be used throughout the paper and some auxiliary results; in Section \ref{sec:interaction} we describe how the diffusions of plane waves interact in the fractional NSE system; in Section \ref{sec:proof} we devote to proving Theorem \ref{Mthm}.

\bigskip

\section{Preliminaries}
\label{sec:pre}

\subsection{Notation}
We denote by $A\lesssim B$ an estimate of the form $A\leq C B$ with
some absolute constant $C$, and by $A\sim B$ an estimate of the form $C_1
B\leq A\leq C_2 B$ with some absolute constants $C_1$, $C_2$.
For simplification of the notation, we denote $\|\cdot\|_{p}=\|\cdot\|_{L^p}$.

\subsection{Semigroup operator $e^{-t(-\Delta)^\alpha}$} 
Consider the Cauchy problem of the $n$ dimensional dissipative equation with a fractional power of the Laplacian,
\begin{equation}\label{eq:FH}
\begin{split}
u_t+(-\Delta)^\alpha u=0, \\
u(x,0)=\phi(x),
\end{split}
\end{equation}
where $(x,t)\in\mathbb{R}^n\times [0,\infty)$.
\noindent Denote by $\mathcal{F}$ and $\mathcal{F}^{-1}$ the Fourier transform and inverse Fourier transform respectively.
Let $e^{-t(-\Delta)^\alpha}$ denote the semigroup generated by a fractional Laplacian:
\begin{equation}\label{sol:FH}
e^{-t(-\Delta)^\alpha}\phi:=\mathcal{F}^{-1}\left(e^{-t|\xi|^{2\alpha}}\mathcal{F}(\phi)(\xi)\right).
\end{equation}
Then $u=e^{-t(-\Delta)^\alpha}\phi$ is a solution of (\ref{eq:FH}).



Let $\mathbb{P}$ denote the projection on divergence-free vector
fields, which acts on a function $\phi$ as
\begin{equation}\notag
\mathbb{P}(\phi)=\phi+\nabla\cdot(-\triangle)^{-1}\div\phi.
\end{equation}
We will use the following well-known estimate.
\begin{Lemma}\label{le:semip}
For any $\phi\in L^\infty$,
\bg\notag
\|\nabla e^{-t(-\Delta)^\alpha}\mathbb P\phi\|_\infty\lesssim t^{-\frac{1}{2\alpha}}\|\phi\|_\infty, \qquad t>0.
\ed

\end{Lemma}

\medskip

\subsection{Norm of Besov spaces}
We recall the definitions of norms for the homogeneous and non-homogeneous Besov spaces 
$\dot{B}_{\infty, \infty}^{-s}$ and $B_{\infty, \infty}^{-s}$ (see \cite{L}) for $s>0$
\begin{align}\label{hombes}
&\|f\|_{\dot{B}_{\infty, \infty}^{-s}}=\sup_{t>0} t^{\frac{s}{2\alpha}}
\|e^{-t(-\Delta)^\alpha}f\|_{L^\infty},\\
&\|f\|_{B_{\infty, \infty}^{-s}} = \sup_{0<t<1}t^{\frac{s}{2\alpha}}
\|e^{-t(-\Delta)^\alpha}f\|_{L^\infty}\notag.
\end{align}

Note that for periodic functions the homogeneous and non-homogeneous norms are equivalent (see \cite{SchT}). 
Therefore, for periodic functions with some fixed period we have
\bg\label{ineq:imbed}
\|f\|_{\dot B_{\infty, \infty}^{-s}} \lesssim \|f\|_{L^\infty},
\ed
since $\|e^{-t(-\Delta)^{\alpha}}f\|_{L^\infty} \leq \|f\|_{L^\infty}$.

We also recall the norm in the Besov space $\dot{B}^{-s}_{\infty,p}$:
\[
\|u\|_{\dot{B}^{-s}_{\infty,p}} = \|\{2^{-sq}\|\Delta_q u\|_\infty \}_{q \in \mathbb{Z}} \|_{l^p},
\]
where $\Delta_q u$ is the Littlewood-Paley projection of $u$.

\medskip

\subsection{Bilinear operator}
Define the bilinear operator
\begin{equation}\label{Ba}
\mathcal{B}_\alpha (u, v)=\int_{0}^{t} e^{-(t-\tau)(-\Delta)^\alpha}
\mathbb{P}\nabla\cdot (u\otimes v)\, d\tau .
\end{equation}
As shown in \cite{KT} in the case of $\alpha=1$, the bilinear operator $\mathcal{B}_1$ continuously
maps $X_T\times X_T$ into $X_T$, where $X_T$ is the Koch-Tataru adapted space. In \cite{BP}, the continuity of $\mathcal B_1$ on $X_T\times X_T$ plays an important role to estimate the higher order iterations (the part $y$ in the paper) of the nonlinear term.
Having to deal with the lack of continuity in the case $\alpha>1$, we use a relatively weaker estimate for the bilinear operator to control
nonlinear interactions (cf. \cite{Yo}):
\begin{Lemma}\label{le:bil} For all $\alpha>0$ the bilinear operator satisfies 
\begin{equation}\label{Bac}
\|\mathcal{B}_\alpha (u, v)\|_{\infty}\lesssim \int_0^t\frac{1}{(t-\tau)^{1/{(2\alpha)}}} \|u(\tau)\|_{\infty}\|v(\tau)\|_{\infty}d\tau.
\end{equation}
\end{Lemma}
\pf
By the definition (\ref{Ba}) and Lemma \ref{le:semip} we have 
\begin{equation}\notag
\begin{split}
\|\mathcal{B}_\alpha (u, v)\|_{\infty}&\lesssim \int_0^t\| e^{-(t-\tau)(-\Delta)^\alpha} \mathbb P\nabla\cdot(u\otimes v)(\tau)\|_\infty d\tau\\
&\lesssim \int_0^t\frac{1}{(t-\tau)^{1/{(2\alpha)}}} \|u(\tau)\|_{\infty}\|v(\tau)\|_{\infty}d\tau.
\end{split}
\end{equation}
\cbdu

\bigskip

\section{Interactions of plane waves}
\label{sec:interaction}

\subsection{The first approximation of a mild solution}
\label{sec:rewrite}
Let $u$ be a solution to (\ref{FNSE}). We write it in the form
\begin{equation}\label{u}
u(t)=e^{-t(-\Delta)^\alpha}u_0-u_1(t)+y(t),
\end{equation}
where
\begin{equation} \label{u1}
u_1(t)= \mathcal{B}_\alpha(e^{-t(-\Delta)^\alpha}u_0,e^{-t(-\Delta)^\alpha}u_0).
\end{equation}
A simple calculation shows that
\begin{equation} \label{y}
y(t) = -\int\limits_{0}^{t}
e^{-(t-\tau) (-\Delta)^\alpha} [G_0 (\tau)+G_1 (\tau)+G_2 (\tau)]d\tau,
\end{equation}
where
\begin{equation}\label{eq:g}
\begin{split}
&G_0=\mathbb{P}[(e^{-t(-\Delta)^\alpha} u_0\cdot\nabla)u_1+(u_1\cdot\nabla)e^{-t(-\Delta)^\alpha} u_0+(u_1\cdot\nabla)u_1]\\
&G_1=\mathbb{P}[(e^{-t(-\Delta)^\alpha} u_0\cdot\nabla)y+(u_1\cdot\nabla)y+(y\cdot\nabla)e^{-t(-\Delta)^\alpha} u_0+(y\cdot\nabla)u_1]\\
&G_2=\mathbb{P}[(y\cdot\nabla)y].
\end{split}
\end{equation}
Note that $G_0$ does not depend on $y$, $G_1$ is linear,
and $G_2$ is quadratic in  $y$.

In this section we show how the diffusions of plane
waves interact in the fractional NSE system.  These interactions are the basis for the
constructions of initial data to produce the norm inflation.
\subsection{Diffusion of a plane wave } As a first step,
we consider the initial data being
one single plane wave. Suppose $k\in \mathbb{R}^3$, $v \in\mathbb{S}^2$ and $k \cdot v=0$. Let
$$
u_0 = v\cos(k\cdot x).
$$
Then
$ \nabla\cdot u_0 = 0$ and
\begin{equation}\label{difu}
e^{-t(-\Delta)^\alpha}v\cos(k\cdot x)=e^{-|k|^{2\alpha}t}v\cos(k\cdot x).
\end{equation}

In fact the ``diffusion"  $e^{-t(-\Delta)^\alpha}v\cos(k\cdot x)$ of a plane wave solves (\ref{FNSE}) with vanishing pressure. And it is important to notice that for $s\geq 0$
$$
\|v\cos(k\cdot x)\|_{\dot{B}^{-s}_{\infty,\infty}} \sim |k|^{-s}.
$$

\subsection{Interaction of plane waves}
\label{sec:int}

Now we consider the interaction of two different single plane waves. Suppose $k_i\in\mathbb{R}^3$, $v_i \in\mathbb{S}^2$ and $k_i\cdot v_i=0$, for $i=1, 2$. Let
\begin{equation}\notag
u_{1}= \cos(k_1\cdot x)v_1, \ \ u_{2}=\cos(k_2\cdot x)v_2.
\end{equation}
To simplify our
calculations we  assume that 
$
k_2\cdot v_1 = \frac 12.
$
It then follows from a straightforward calculation that
\begin{equation}\notag
\begin{split}
&e^{-t(- \Delta)^\alpha} u_{1}\cdot\nabla(e^{-t (-\Delta)^\alpha} u_{2})\\
&=- e^{-(|k_1|^{2\alpha}+|k_2|^{2\alpha}) t}v_2\cos(k_1\cdot x)\sin(k_2\cdot x)(k_2\cdot v_1)\\
&=-\frac{1}{4} e^{-(|k_1|^{2\alpha}+|k_2|^{2\alpha})t}v_1 (\sin((k_2-k_1)\cdot x) +
\sin((k_1+k_2)\cdot x)). 
\end{split}
\end{equation}
Hence
\begin{equation}\notag
\begin{split}
&\mathcal B_\alpha(e^{-t(- \Delta)^\alpha} u_{1}, e^{-t(- \Delta)^\alpha} u_{2}) \\
=& \frac 14  v_1\sin((k_2-k_1)\cdot x) \int_0^t
e^{-(|k_1|^{2\alpha}+|k_2|^{2\alpha})\tau}e^{-|k_2-k_1|^{2\alpha}(t-\tau)}d\tau \\
 & + \frac
14 v_1 \sin((k_1+k_2)\cdot x)\int_0^t
e^{-(|k_1|^{2\alpha}+|k_2|^{2\alpha})\tau}e^{-|k_1+k_2|^{2\alpha}(t-\tau)}d\tau.
\end{split}
\end{equation}
Therefore,
the interaction of the two plane waves is small 
in $\dot{B}^{-s}_{\infty, \infty}$ 
if neither the sum nor the difference
of their wave vectors is small in magnitude. In the contrary, the
interaction is sizable in $\dot{B}^{-s}_{\infty, \infty}$ if either the
sum or the difference of their wave vectors is small in magnitude.

\bigskip

\section{Proof of theorem \ref{Mthm}}
\label{sec:proof}
In this section we follow the idea from \cite{BP} to construct 
initial data that produce norm inflation for solutions to the fractional Navier-Stokes equation.
From the discussions in Subsection \ref{sec:int} it is clear that the
interaction of two plane waves is not enough to produce the norm
inflation, which actually requires a large number of waves.
We also make sure that the initial data is space-periodic and smooth,
which ensures the local existence of a smooth periodic solution to the
fractional NSE. As we control its $L^{\infty}$ norm, the solution will
remain smooth until the time of the norm inflation.

\subsection{Construction of initial data for the fractional NSE system}

For a fixed small number $\delta>0$, 
the initial data will be chosen as follows:
\begin{equation}\label{u0}
u_0=r^{-\beta}\sum_{i=1}^r|k_i|^{\alpha}(v \cos(k_i\cdot x)+v' \cos(k'_i\cdot x)),
\end{equation}
where $\beta>0$. We expect for each $i$ the interaction of the two plane waves
$v\cos(k_i\cdot x)$ and $v'\cos(k_i'\cdot x)$ to be sizable in
$\dot{B}^{-s}_{\infty, \infty}$, while the interactions of plane waves corresponding to different
indexes $i$ to be small.
 Hence, we choose
\begin{itemize}
\item{Wave vectors:} Let $\zeta=(1,0,0)$ and $\eta=(0,0,1)$. The wave vectors $k_i\in\mathbb{Z}^3$ are parallel to $\zeta$. Let $K$ be a large integer dependent on $r$. The magnitude of $k_i$ is defined by
\begin{align}\label{kl}
|k_i|=2^{i-1}K, \ \ \ \ i=1, 2, 3, ..., r.
\end{align}
The wave vectors $k'_i\in\mathbb{Z}^3$ are defined by
\begin{equation} \label{ks-ks'}
k'_i=k_i+\eta.
\end{equation}
\item{Amplitude vectors:} Let 
\bg\label{eq:v}
v=(0,0,1), \ \ \ v'=(0,1,0). 
\ed
Hence
\begin{equation}\notag
k_i\cdot v=k'_i\cdot v'=0,
\end{equation}
which ensures that the initial data is divergence free.
\end{itemize}

We first point out  the following  simple
facts to further motivate the choices of the parameters.

\begin{Lemma}\label{k}Let $\gamma> 0$, $\alpha\geq 1$. With the choices (\ref{kl})-(\ref{eq:v}), the following holds:
\begin{equation}\label{vk0}
k_i \cdot v' = 0,\ \ \ \ k'_i\cdot v =1, \quad \forall
\quad i = 1,2, \dots, r,
\end{equation}
\begin{equation}\label{kll1}
\sum_{j<i}|k_j|^{\alpha}\sim |k_{i-1}|^{\alpha} \quad\text{and } \ \sum_{j<i}|k_j'|^{\alpha}\sim
|k_{i-1}'|^{\alpha},
\end{equation}
\begin{equation}\label{ksum2}
\sum_{i=1}^r |k_i|^{\gamma}e^{-|k_i|^{2\alpha}t} \lesssim t^{-\frac{\gamma}{2\alpha}}\quad
\text{and } \ \sum_{i=1}^r |k_i'|^{\gamma}e^{-|k_i'|^{2\alpha}t} \lesssim t^{-\frac{\gamma}{2\alpha}} .
\end{equation}
\end{Lemma}
\pf The first conclusion (\ref{vk0}) is obvious due to (\ref{kl})-(\ref{eq:v}). 
By the definition (\ref{kl}), it is clear that
$|k_{i-1}|^{\alpha}<\frac{1}{2}|k_i|^{\alpha}$, which immediately implies (\ref{kll1}). 
Thanks to (\ref{kl}), we have that $|k_i|^{\alpha}\sim |k_i|^{\alpha}-|k_{i-1}|^{\alpha}$. Thus,
$$
\sum_{i=1}^r|k_i|^{\gamma}e^{-|k_i|^{2\alpha}t} \sim
\sum_{i=1}^r|k_i|^{\gamma-\alpha}(|k_i|^{\alpha}-|k_{i-1}|^{\alpha})e^{-|k_i|^{2\alpha}t},
$$
while the latter one can be considered as a finite Riemman summation
of the function $x^{\gamma/\alpha-1}e^{-x^2t}$. Therefore, for $\gamma>0$ and $\alpha>0$,
\begin{equation}\notag
\sum_{i=1}^r|k_i|^{\gamma}e^{-|k_i|^{2\alpha}t} \lesssim \int_0^\infty x^{\gamma/\alpha-1}e^{-x^2t} dx
=t^{-\frac{\gamma}{2\alpha}}\int_0^\infty y^{\gamma/\alpha-1}e^{-y^2} dy \lesssim t^{-\frac{\gamma}{2\alpha}}.
\end{equation}
\cbdu

Next we estimate the norms of the initial data.
\begin{Lemma} \label{le:u0}
Let $u_0$ be given in (\ref{u0}) and $\alpha>0$. Then
\begin{equation}\label{norm:u0-1}
\|u_0\|_{\dot{B}_{\infty, p}^{-\alpha}}\lesssim r^{1/p-\beta}, \qquad 1\leq p \leq \infty.
\end{equation}
\end{Lemma}
\pf Due to (\ref{difu}), we have that, 
\begin{equation} \label{h-u0}
e^{-t(-\Delta)^\alpha}u_0=r^{-\beta}\sum_{i=1}^r|k_i|^{\alpha}(v
\cos(k_i\cdot x)e^{-|k_i|^{2\alpha}t}+v'\cos(k'_i\cdot x)e^{-|k'_i|^{2\alpha}t}).
\end{equation}
Hence by Lemma \ref{k}, 
\bg \label{first-bound}
\|u_0\|_{\dot{B}_{\infty, \infty}^{-\alpha}} \sim r^{-\beta}
\sup_{0<t<1} t^{\frac{1}{2}}\sum_{i=1}^r |k_i|^{\alpha}\left(e^{-|k_i|^{2\alpha}t}+e^{-|k'_i|^{2\alpha}t}\right)
\lesssim r^{-\beta}\notag.
\ed
A direct computation also gives
\[
\|u_0\|_{\dot{B}_{\infty, p}^{-\alpha}} \lesssim r^{-\beta} \left( \sum_{i=1}^r 1^{p} \right)^{1/p} = r^{1/p - \beta}, \qquad p \geq1.
\]
 \cbdu
\begin{Lemma}\label{le:u0infty} Let $u_0$ be given in (\ref{u0}). Then
\begin{equation}\notag
\|e^{-t(-\Delta)^\alpha}u_0\|_{\infty} \lesssim r^{-\beta}t^{-1/2}.
\end{equation}
\end{Lemma}
\pf By (\ref{h-u0}) and Lemma \ref{k}, we infer that
\bg\notag
\|e^{-t(-\Delta)^\alpha}u_0\|_{\infty} \lesssim r^{-\beta}
 \sum_{i=1}^r |k_i|^{\alpha}\left(e^{-|k_i|^{2\alpha}t}+e^{-|k'_i|^{2\alpha}t}\right)
 \lesssim r^{-\beta}t^{-1/2}.
\ed
\cbdu

\bigskip

\subsection{Analysis of $u_1$}

As demonstrated in Subsection \ref{sec:rewrite} we consider the decomposition
$$
\aligned u  & = e^{-t(-\Delta)^\alpha}u_0 - u_1 + y. \endaligned
$$
Recall the definition (\ref{u1})
$$
u_1 = \mathcal{B}_\alpha(e^{-t(-\Delta)^\alpha}u_0, e^{-t(-\Delta)^\alpha}u_0). 
$$
By (\ref{eq:v}), (\ref{vk0}), (\ref{h-u0}) and a straightforward calculation, it follows that
\begin{equation}\label{u0interaction}
\begin{split}
&(e^{-t(-\Delta)^\alpha}u_0\cdot\nabla) e^{-t(-\Delta)^\alpha}u_0 \\
 =&-r^{-2\beta}\sum_{i=1}^r\sum_{j=1}^r|k_i|^{\alpha}|k_j|^{\alpha}e^{-(|k_i|^{2\alpha}+|k'_j|^{2\alpha})t}v'\cos(k_i\cdot x)\sin(k'_j\cdot x)\\ 
=&-\frac{r^{-2\beta}}{2}\sum_{i=1}^r|k_i|^{2\alpha}e^{-(|k_i'|^{2\alpha}+|k_i|^{2\alpha})t}\sin(\eta\cdot x)v'\\ 
& -\frac{r^{-2\beta}}{2}\sum_{i\neq j}^r|k_i|^{\alpha}|k_j|^{\alpha}e^{-(|k_i|^{2\alpha}+|k'_j|^{2\alpha})t}\sin((k'_j-k_i)\cdot x)v'\\ 
& -\frac{r^{-2\beta}}{2}\sum_{i=1}^r\sum_{j=1}^r|k_i|^{\alpha}|k_j|^{\alpha}e^{-(|k_i|^{2\alpha}+|k'_j|^{2\alpha})t}\sin((k'_j+k_i)\cdot x)v'\\
\equiv &E_0+E_1+E_2,
\end{split}
\end{equation}
where we used the formula $\cos x\sin y=[\sin(x+y)-\sin(x-y)]/2$.\\

Recall that $\eta\cdot v'=0$, $(k'_j+k_i)\cdot v'=0$ and $(k'_j-k_i)\cdot v'=0$ for all $i,j$ due to (\ref{vk0}). Hence $E_{0}$, $E_{1}$ and $E_{2}$ are divergence free vectors. 
Thus we can write 
\begin{equation}\label{dec:u1}
\begin{split}
u_1=&\int_0^te^{-(t-\tau)(-\Delta)^\alpha}E_{0}(\tau)d\tau+\int_0^te^{-(t-\tau)(-\Delta)^\alpha}E_{1}(\tau)d\tau\\
&+\int_0^te^{-(t-\tau)(-\Delta)^\alpha}E_{2}(\tau)d\tau\equiv u_{10}+u_{11}+u_{12}.
\end{split}
\end{equation}

We have the following estimates.

\begin{Lemma}\label{le:u10}
Let $u_{10}$ be defined in (\ref{dec:u1}) and $s> 0$. Then
\begin{align}\notag 
&\|u_{10}(\cdot, t)\|_{\dot B_{\infty,\infty}^{-s}}\gtrsim r^{1-2\beta}, \qquad \mbox {for all} \qquad   K^{-2\alpha}\leq t\leq 1,\\
&\|u_{10}(\cdot, t)\|_{\infty} \lesssim r^{1-2\beta}, \qquad \mbox {for all} \qquad t>0\notag.
\end{align}
\end{Lemma}
\pf
From (\ref{u0interaction}) and (\ref{dec:u1}) it follows by a straightforward calculation
\begin{equation}\notag
\begin{split}
u_{10}&=-\frac{r^{-2\beta}}{2}\int_0^t\sum_{i=1}^r|k_i|^{2\alpha}e^{-(|k_i'|^{2\alpha}+|k_i|^{2\alpha})\tau}e^{-|\eta|^{2\alpha}(t-\tau)}\sin(\eta\cdot x)v'd\tau\\
&=-\frac{r^{-2\beta}}{2}\sin (\eta\cdot x)v'\sum_{i=1}^r|k_i|^{2\alpha}e^{-t}\frac{1-e^{-(|k_i'|^{2\alpha}+|k_i|^{2\alpha}-1)t}}{|k'_i|^{2\alpha}+|k_i|^{2\alpha}-1}\\
&\sim -\frac{r^{-2\beta}}{2}\sin (\eta\cdot x)v'\sum_{i=1}^re^{-t}(1-e^{-|k_i|^{2\alpha}t}).
\end{split}
\end{equation}
 Hence for $K^{-2\alpha}\leq t\leq 1$ and $s> 0$,
\bg\notag
\|u_{10}(\cdot,t)\|_{\dot B_{\infty,\infty}^{-s}}\gtrsim r^{-2\beta}\cdot r\sup_{0<\tau<1}\tau^{\frac{s}{2\alpha}}e^{-|\eta|^{2\alpha}\tau}
\gtrsim r^{1-2\beta}.
\ed
On the other hand,
\bg\notag
\|u_{10}(\cdot,t)\|_{\infty}\lesssim \frac{r^{-2\beta}}{2}\cdot r\lesssim r^{1-2\beta},
\ed
for all $t>0$.
\cbdu

\begin{Lemma}\label{le:u11}
Let $u_{11}$ and $u_{12}$ be defined in (\ref{dec:u1}). Then
\begin{equation}\notag
\|u_{11}(\cdot, t)\|_{\infty}\lesssim r^{-2\beta}, \quad
\|u_{12}(\cdot, t)\|_{\infty} \lesssim r^{-2\beta},
\end{equation}
for all $t>0$.
\end{Lemma}
\pf
Thanks to (\ref{u0interaction}) and (\ref{dec:u1}), it follows that
\begin{equation}\notag
\begin{split}
u_{11}=&\frac{r^{-2\beta}}{2}\int_0^t\sum_{i\neq j}^r|k_i|^{\alpha}|k_j|^{\alpha}e^{-(|k_i|^{2\alpha}+|k'_j|^{2\alpha})\tau}e^{-|k'_j-k_i|^{2\alpha}(t-\tau)}
\sin((k'_j-k_i)\cdot x)v'd\tau\\
\sim& \frac{r^{-2\beta}}{2}\sum_{i=1}^r\sum_{j<i}|k_i|^{\alpha}|k_j|^{\alpha}e^{-|k_i-k'_j|^{2\alpha}t}\frac{1-e^{-(|k_i|^{2\alpha}+|k'_j|^{2\alpha}-|k_i-k'_j|^{2\alpha})t}}{|k_i|^{2\alpha}+|k'_j|^{2\alpha}-|k_i-k'_j|^{2\alpha}}\\
&\cdot\sin ((k'_j-k_i)\cdot x)v'\\
\sim& \frac{r^{-2\beta}}{2}\sum_{i=1}^r\sum_{j<i}|k_i|^{\alpha}|k_j|^{\alpha}te^{-|k_i|^{2\alpha}t}
\sin ((k'_j-k_i)\cdot x)v',
\end{split}
\end{equation}
where we used the fact that $\frac{1-e^{-x}}{x}$ is bounded for $x>0$. 
Hence, by (\ref{kll1}) and (\ref{ksum2}) we infer that
\begin{equation}\notag
\begin{split}
\|u_{11}(\cdot, t)\|_{\infty}&\lesssim r^{-2\beta}\sum_{i=1}^r\sum_{j<i}|k_i|^{\alpha}|k_j|^{\alpha}te^{-|k_i|^{2\alpha}t}\\
&\lesssim r^{-2\beta}\sum_{i=1}^r|k_i|^{2\alpha}te^{-|k_i|^{2\alpha}t}\lesssim r^{-2\beta}.
\end{split}
\end{equation}
Similarly, we have
\begin{equation}\notag
\begin{split}
u_{12}=&\frac{r^{-2\beta}}{2}\int_0^t\sum_{i=1}^r\sum_{j=1}^r|k_i|^{\alpha}|k_j|^{\alpha}e^{-(|k_i|^{2\alpha}+|k'_j|^{2\alpha})\tau}e^{-|k_i+k'_j|^{2\alpha}(t-\tau)}\sin((k_i+k'_j)\cdot x)v'd\tau\\
=&\frac{r^{-2\beta}}{2}
\sum_{i=1}^r\sum_{j=1}^r|k_i|^{\alpha}|k_j|^{\alpha}e^{-(|k_i|^{2\alpha}+|k'_j|^{2\alpha})t}\frac{1-e^{-(|k_i+k'_j|^{2\alpha}-|k_i|^{2\alpha}-|k'_j|^{2\alpha})t}}{|k_i+k'_j|^{2\alpha}-|k_i|^{2\alpha}-|k'_j|^{2\alpha}}\\
&\cdot\sin ((k_i+k'_j)\cdot x)v'\\
\sim& r^{-2\beta}\sum_{i=1}^r\sum_{j\leq i}|k_i|^{\alpha}|k_j|^{\alpha}e^{-(|k_i|^{2\alpha}+|k'_j|^{2\alpha})t}t\sin ((k_i+k'_j)\cdot x)v'.
\end{split}
\end{equation}
Thus,
\begin{equation}\notag
\begin{split}
\|u_{12}(\cdot, t)\|_{\infty}&\lesssim r^{-2\beta}\sum_{i=1}^r\sum_{j\leq i}|k_i|^{\alpha}|k_j|^{\alpha}te^{-|k_i|^{2\alpha}t}\\
&\lesssim r^{-2\beta}\sum_{i=1}^r|k_i|^{2\alpha}te^{-|k_i|^{2\alpha}t}\lesssim r^{-2\beta}.
\end{split}
\end{equation}
\cbdu

\bigskip

\subsection{Analysis of $y$}

In this section we analyze the part $y$ of the solution.
The idea is to control $y$ using the estimate (\ref{Bac}) of the
bilinear operator $\mathcal B_\alpha$ in the space $L^\infty$. \\

Recall from Subsection \ref{sec:rewrite} that
\bg\label{y-t}
y(t) = -\int\limits_{0}^{t}
e^{-(t-\tau) (-\Delta)^\alpha} [G_0 (\tau)+G_1 (\tau)+G_2 (\tau)]d\tau. 
\ed 

\begin{Lemma}\label{le:yinfty}
Let $\alpha\geq 1$ and $\beta \in (0,1/2)$. Then
\bg\notag
\|y(t)\|_{\infty}\lesssim r^{1-3\beta}t^{\frac{1}{2}-\frac{1}{2\alpha}}+r^{2-4\beta}t^{1-\frac{1}{2\alpha}},
\qquad \forall t \in[0,T],
\ed
provided  $T$ is small and $r$ is large enough.
\end{Lemma}
\pf 
It follows from (\ref{eq:g}) and (\ref{y-t}) that
\begin{equation}\notag
\begin{split}
\|y(t)\|_\infty&\lesssim \|\mathcal B_\alpha(e^{-t(-\Delta)^\alpha}u_0, u_1)\|_\infty+
\|\mathcal B_\alpha(u_1, u_1)\|_\infty\\
&+\|\mathcal B_\alpha(e^{-t(-\Delta)^\alpha}u_0, y)\|_\infty+
\|\mathcal B_\alpha(u_1, y)\|_\infty+\|\mathcal B_\alpha(y, y)\|_\infty.
\end{split}
\end{equation}

Applying the bilinear estimate (\ref{Bac}), Lemmas \ref{le:u0infty}, \ref{le:u10}, and \ref{le:u11} we infer
\begin{equation}\notag
\begin{split}
\|\mathcal B_\alpha(e^{-t(-\Delta)^\alpha}u_0, u_1)\|_\infty
&\lesssim\int_0^t\frac{1}{(t-\tau)^{1/{(2\alpha)}}}\|e^{-\tau(-\Delta)^\alpha}u_0\|_{\infty}\|u_1(\tau)\|_{\infty}d\tau\\
&\lesssim r^{1-3\beta}\int_0^t(t-\tau)^{-1/{(2\alpha)}}\tau^{-1/2}d\tau\\
&\lesssim r^{1-3\beta}t^{\frac{1}{2}-\frac{1}{2\alpha}},
\end{split}
\end{equation}
where we used the boundedness of Beta function for $\alpha>1/2$:
$$
\int_0^t(t-\tau)^{-1/{(2\alpha)}}\tau^{-1/2}d\tau=t^{\frac{1}{2}-\frac{1}{2\alpha}}B(\frac{1}{2}, 1-\frac{1}{2\alpha})\leq Ct^{\frac{1}{2}-\frac{1}{2\alpha}}.
$$
Similarly, using the estimates obtained in previous two subsections, we obtain
\begin{equation}\notag
\|\mathcal B_\alpha(u_1, u_1)\|_\infty
\lesssim\int_0^t\frac{1}{(t-\tau)^{1/{(2\alpha)}}}\|u_1(\tau)\|_{\infty}^2d\tau\lesssim r^{2-4\beta}
t^{1-\frac{1}{2\alpha}},
\end{equation}
\begin{equation}\notag
\begin{split}
\|\mathcal B_\alpha(e^{-t(-\Delta)^\alpha}u_0, y)\|_\infty
&\lesssim\int_0^t\frac{1}{(t-\tau)^{1/{(2\alpha)}}}\|e^{-\tau(-\Delta)^\alpha}u_0\|_{\infty}\|y(\tau)\|_{\infty}d\tau\\
&\lesssim r^{-\beta}\int_0^t(t-\tau)^{-1/{(2\alpha)}}\tau^{-1/2}d\tau\sup_{0<\tau<t}\|y(\tau)\|_{\infty}  \\
&\lesssim r^{-\beta}t^{\frac{1}{2}-\frac{1}{2\alpha}}\sup_{0<\tau<t}\|y(\tau)\|_{\infty},
\end{split}
\end{equation}
\begin{equation}\notag
\begin{split}
\|\mathcal B_\alpha(u_1, y)\|_\infty
&\lesssim\int_0^t\frac{1}{(t-\tau)^{1/{(2\alpha)}}}\|u_1(\tau)\|_{\infty}\|y(\tau)\|_{\infty}d\tau\\
&\lesssim r^{1-2\beta}\int_0^t(t-\tau)^{-1/{(2\alpha)}}d\tau\sup_{0<\tau<t}\|y(\tau)\|_{\infty}\\
&\lesssim r^{1-2\beta}t^{1-\frac{1}{2\alpha}}\sup_{0<\tau<t}\|y(\tau)\|_{\infty},
\end{split}
\end{equation}
\begin{equation}\notag
\|\mathcal B_\alpha(y, y)\|_\infty
\lesssim\int_0^t\frac{1}{(t-\tau)^{1/{(2\alpha)}}}\|y(\tau)\|_{\infty}^2d\tau
\lesssim t^{1-\frac{1}{2\alpha}}\left(\sup_{0<\tau<t}\|y(\tau)\|_{\infty}\right)^2\notag.
\end{equation}
Thus we have
\begin{align}\notag
\begin{split}
\|y(t)\|_{\infty}
&\lesssim r^{1-3\beta}t^{\frac{1}{2}-\frac{1}{2\alpha}}+r^{2-4\beta}t^{1-\frac{1}{2\alpha}}\\
&+\left(r^{-\beta}t^{\frac{1}{2}-\frac{1}{2\alpha}}+r^{1-2\beta}t^{1-\frac{1}{2\alpha}}+
t^{1-\frac{1}{2\alpha}}\sup_{0<\tau<t}\|y(\tau)\|_{\infty}\right)\sup_{0<\tau<t}\|y(\tau)\|_{\infty}.
\end{split}
\end{align}

We choose large enough $r$ and  small enough $T>0$,  such that 
\begin{equation}\label{para1}
A:=r^{-\beta}t^{\frac{1}{2}-\frac{1}{2\alpha}}+r^{1-2\beta}t^{1-\frac{1}{2\alpha}}+
t^{1-\frac{1}{2\alpha}}(r^{1-3\beta}t^{\frac{1}{2}-\frac{1}{2\alpha}}+r^{2-4\beta}t^{1-\frac{1}{2\alpha}})\ll 1
\end{equation}
for $0\leq t\leq T$. Indeed, note that the powers of $t$ in $A$ are all nonnegative for $\alpha\geq 1$. Thus,
\begin{equation}\notag
A\leq r^{-\beta}T^{\frac{1}{2}-\frac{1}{2\alpha}}+r^{1-2\beta}T^{1-\frac{1}{2\alpha}}+
T^{1-\frac{1}{2\alpha}}(r^{1-3\beta}T^{\frac{1}{2}-\frac{1}{2\alpha}}+r^{2-4\beta}T^{1-\frac{1}{2\alpha}}).
\end{equation}
 Let $T=r^{-\gamma}$. It follows 
 \begin{equation}\label{para-A}
A\leq r^{-\beta-\gamma(\frac{1}{2}-\frac{1}{2\alpha})}+r^{1-2\beta-\gamma(1-\frac{1}{2\alpha})}+
r^{1-3\beta-\gamma(\frac{3}{2}-\frac{1}{\alpha})}+r^{2-4\beta-\gamma(2-\frac{1}{\alpha})}.
\end{equation}
We choose $\gamma$ such that 
\bg\label{para-gamma}
\gamma>\frac{1-2\beta}{1-1/(2\alpha)},
\ed
which guarantees all the powers of $r$ are negative in (\ref{para-A}). Hence (\ref{para1}) is satisfied for $r$ large enough. Since $y(0)=0$, we have the following bound by an absorbing argument:
\bg\notag
\|y(t)\|_\infty\lesssim r^{1-3\beta}t^{\frac{1}{2}-\frac{1}{2\alpha}}+r^{2-4\beta}t^{1-\frac{1}{2\alpha}},
\ed
for all $0<t\leq T$. 
\cbdu

\subsection{Finishing the proof}
\label{subsec:end}

Now we are ready to complete the proof of Theorem \ref{Mthm}.
Since $u_0$ is smooth and space-periodic, there exists $T^*>0$ and
a smooth space-periodic solution $u(t)$ to \eqref{FNSE} on $[0,T^*)$
with $u(0)=u_0$, such that either $T^*= +\infty$ or
\[
\limsup_{t\to T^*-} \|u(t)\|_{\infty} = +\infty.
\]
Lemmas~\ref{le:u10},
\ref{le:u11}, and \ref{le:yinfty} imply that $T^* >T$, where $T=r^{-\gamma}$ and $\gamma$ is large enough so that \eqref{para-gamma}
holds. Note that $T<1$.

Now using
(\ref{u}), we combine the imbedding estimate (\ref{ineq:imbed}), Lemmas \ref{le:u0infty}, \ref{le:u10}, \ref{le:u11} and \ref{le:yinfty} to obtain that, for $K^{-2\alpha}\leq t\leq T$,
\begin{equation}\label{ineq:blow-up}
\begin{split}
\|u(\cdot, t)\|_{\dot{B}_{\infty,\infty}^{-s}}
\geq&\|u_{10}(\cdot, t)\|_{\dot{B}_{\infty,\infty}^{-s}} - \|u_{11}(\cdot, t)\|_{\infty}-\|u_{12}(\cdot, t)\|_{\infty}\\
&-\|e^{-t(-\Delta)^{\alpha}}u_0\|_{\infty}
- \|y(\cdot, t)\|_{\infty}  \\
 \gtrsim &r^{1-2\beta}\left(1-r^{-1}-r^{\beta-1}t^{-\frac{1}{2}}-r^{-\beta}t^{\frac{1}{2}-\frac{1}{2\alpha}}-r^{1-2\beta}t^{1-\frac{1}{2\alpha}}\right)\\
 \gtrsim &r^{1-2\beta}\left(1-r^{\beta-1}K^{\alpha}-r^{-\beta}T^{\frac{1}{2}-\frac{1}{2\alpha}}-r^{1-2\beta}T^{1-\frac{1}{2\alpha}}\right).
\end{split}
\end{equation}
We will show that we can choose parameters so that 
\begin{equation}\label{small-p}
B:=r^{\beta-1}K^{\alpha}+r^{-\beta}T^{\frac{1}{2}-\frac{1}{2\alpha}}+r^{1-2\beta}T^{1-\frac{1}{2\alpha}}\leq1/4, \qquad \mbox {for} \quad \alpha\geq 1.
\end{equation}
Let $K=r^{\zeta}$ with positive $\zeta$, and recall that $T=r^{-\gamma}$ as in Lemma \ref{le:yinfty}. Then
\[
B=r^{\beta-1+\zeta\alpha}+r^{-\beta-\gamma(\frac{1}{2}-\frac{1}{2\alpha})}+r^{1-2\beta-\gamma(1-\frac{1}{2\alpha})}.
\]
For any $\beta\in(0,\frac{1}{2})$ we choose $\zeta,\gamma$ such that 
\bg\label{para-zeta}
0<\zeta<\frac{1-\beta}{\alpha}, \qquad \frac{1-2\beta}{1-1/(2\alpha)} < \gamma <2\alpha\zeta,
\ed
which can be done because $\alpha \geq 1$. 
This implies that all the powers of $r$ in $B$ are negative and hence (\ref{small-p}) is satisfied for $r$ large enough.
Moreover, the condition $\gamma<2\alpha\zeta$ guarantees that  $K^{-2\alpha}< T$. 
Note that the conditions on $\gamma$ in (\ref{para-gamma}) and (\ref{para-zeta}) coincide. 


Given any $\delta>0$ in Theorem \ref{Mthm}, we now choose a suitable large $r$ such that
\bg\notag
r^{1-2\beta}\gtrsim \frac{1}{\delta}. 
\ed
Therefore, it follows from (\ref{ineq:blow-up}) and (\ref{small-p}) that
\begin{equation}\notag
\|u(T)\|_{\dot{B}_{\infty,\infty}^{-s}} \gtrsim r^{1-2\beta}\gtrsim \frac{1}{\delta}.
\end{equation}
Finally, Lemma \ref{le:u0} implies
that the initial data $u_0$ satisfies
\bg\notag
\|u_0\|_{\dot{B}_{\infty,p}^{-\alpha}}\lesssim r^{1/p-\beta} \lesssim r^{2\beta-1} \lesssim \delta, 
\ed
as long as $3\beta \geq 1+1/p$, which holds for any $p>2$ provided $\beta$ is close enough to $1/2$. 
This competes the proof of Theorem \ref{Mthm}.


\bigskip


\end{document}